\documentclass[8pt]{article}
\usepackage[francais,english]{babel}
\usepackage[T1]{fontenc}
\usepackage{kpfonts}
\usepackage{color}\usepackage{color}
\usepackage{bera}
\usepackage{amssymb}
\usepackage{geometry}
\usepackage[utf8]{inputenc}
\usepackage{amsmath,mathtools}
\usepackage{graphicx}
\usepackage{sectsty}
\usepackage{lipsum}
\usepackage{mathrsfs}
\usepackage{amsthm}
\usepackage{multicol}
\usepackage{comment}
\usepackage{wrapfig}
\usepackage{float}
\usepackage{lipsum}

\usepackage{listings}

\usepackage[numbers,sort&compress]{natbib}
\allsectionsfont{\centering}
\newtheorem{theorem}{Theorem}[section]

\newtheorem{example}{Example}[section]
\newtheorem{definition}{Definition}[section]

\newtheorem{lemma}{Lemma}[section]

\numberwithin{equation}{section}
\begin{document}
\large

\begin{center}
	{\fontsize{16}{24}\bf A Predictor-Corrector Algorithm in the Framework of Conformable Fractional Differential Equations
		
		%
		%
	}\end{center}
\vspace{.3cm}
   \hrule
   \begin{center}
   \textbf{Mohamed Echchehira }$^{(1,a)}$,
   \textbf{Youness Assebbane }$^{(1,b)}$, \textbf{Mustapha Atraoui}$^{(1,c)}$ and \textbf{Mohamed Bouaouid}$^{(2,d)}$
\end{center}
\begin{center}
$^{(1)}$ Ibn Zohr University, Ait Melloul University Campus, Faculty of Applied Sciences, Department of
Mathematics, Research Laboratory for Innovation in Mathematics and Intelligent Systems, BP 6146, 86150, Agadir, Morocco.\\
$^{(2)}$Sultan Moulay Slimane University, National School of Applied Sciences, 23000 Béni-Mellal, Morocco.
\end{center}
\begin{center}
$^{(a)}$ mohamed.echchehira.57@edu.uiz.ac.ma
\end{center}
\begin{center}
$^{(b)}$ youness.assebbane.97@edu.uiz.ac.ma
\end{center}
\begin{center}
$^{(c)}$ m.atraoui@uiz.ac.ma
\end{center}
\begin{center}
$^{(d)}$ bouaouidfst@gmail.com
\end{center}

\begin{center}
	{\fontsize{14}{24}\bf Abstract
		
		%
		%
	}\end{center}
This work proposes a conformable fractional predictor-corrector algorithm for solving conformable fractional differential equations. Fractional calculus is finding applications in various scientific fields, but existing numerical methods might have limitations. This work addresses that gap by introducing a new algorithm specifically designed for the conformable fractional derivative using Adams-Bashforth and Adams-Moulton methods.
\begin{description}
  \item 
  \item[Keywords:] Fractional differential equations; Conformable fractional derivative; Adams-Bashforth and Adams-Moulton methods; Numerical solution.
\end{description}

\section{Introduction}
The theory of fractional derivatives is an ancient concept, originating from a conversation on September 30, 1695, between L'Hôpital and Leibniz. They discussed the definition of the operator $\frac{d^n(.)}{dt^n}$ for $n=\frac{1}{2}$, pondering "what if $n$ is fractional?" \cite{butzer2000introduction}. Despite its long history, fractional calculus is both an old and novel topic. It has only been the subject of specialized conferences and treatises since the 1970s. Shortly after completing his Ph.D. dissertation on fractional calculus, B. Ross organized the First Conference on Fractional Calculus and its Applications at the University of New Haven in June 1974 and edited its proceedings \cite{machado2011recent}.
\\\\
Fractional calculus, or the calculus of non-integer order differentiation and integration, has garnered significant interest due to its ability to more accurately model and solve complex phenomena in various scientific fields \cite{shaikh2020mathematical,magin2010fractional,dalir2010applications,koh1990application}. For instance, in speech signal modeling, fractional calculus provides a superior alternative to the traditional Linear Predictive Coding (LPC) approach by employing fractional order integrals as basis functions, thus requiring fewer parameters and achieving more accurate signal representation \cite{assaleh2007modeling}. One of the most extensive application areas of fractional calculus has been linear viscoelasticity, due to its effectiveness in modeling hereditary phenomena with long memory, for more details, refer to the book \cite{mainardi2022fractional} and to the works \cite{bagley1983theoretical,meral2010fractional,koeller1984applications}. Similarly, in the realm of acoustics, fractional derivatives have been successfully applied to describe sound wave propagation through porous materials, offering a new predictive method validated through experimental work on various media such as plastic foams and sandy sediments \cite{fellah2002application}. Furthermore, in image processing, edge detection can be significantly enhanced by fractional differentiation, which improves detection selectivity and robustness to noise, outperforming traditional integer-order methods \cite{mathieu2003fractional}. Additionally, fluid mechanics has benefited from fractional calculus, where its application to time-dependent viscous-diffusion problems has yielded simpler and more powerful solutions compared to classical methods, thereby validating its effectiveness and broad applicability \cite{kulish2002application}.
\\\\
Over time, multiple approaches have been presented in the literature to define the concept of fractional derivative, such as the Riemann-Liouville definition and the Caputo definition. For $n-1 \leq \alpha < n$ with $n\in \mathbb{Z}^{*}_{+}$, the fractional derivative of a function $f$ in the Riemann-Liouville sense is given by \cite{podlubny1998fractional}
\begin{equation}
    ^{RL}D_a^\alpha f(t)=\frac{1}{\Gamma(n-\alpha)} \frac{d^n}{d t^n} \int_a^t \frac{f(x)}{(t-x)^{\alpha-n+1}} d x,
\end{equation}
where $\Gamma$ is the Euler Gamma function. We note that the Riemann-Liouville derivative of a constant function is not necessarily zero, but this property holds true for the fractional derivative in the Caputo sense, defined by \cite{podlubny1998fractional}
\begin{equation}
    ^{C}D_a^\alpha f(t)=\frac{1}{\Gamma(n-\alpha)} \int_a^t \frac{f^{(n)}(x)}{(t-x)^{\alpha-n+1}} d x.\label{caputoOperator}
\end{equation}
In their work, Khalil et al. introduce a novel definition of the fractional derivative and integral, presenting it as the most natural and fruitful formulation to date. This new definition for $0 \leq \alpha < 1$ aligns with classical definitions when applied to polynomials, differing only by a constant. Moreover, when $\alpha = 1$, it seamlessly matches the classical definition of the first derivative. This innovative approach not only offers theoretical elegance but also proves to be highly practical, with the authors demonstrating its applicability to various fractional differential equations \cite{khalil2014new}.
\begin{definition}Let $\alpha\in (0,1]$, the conformable fractional derivative of order $\alpha$
of a function $f:\hspace{0.2cm}[0,+\infty)\longrightarrow\mathbb{R}$ is defined as
\begin{equation}
T_{\alpha}f(t)=\displaystyle{\lim_{\varepsilon\longrightarrow0}\frac{f(t+\varepsilon t^{1-\alpha})-f(t)}{\varepsilon}}\hspace{0.2cm} \mbox{for} \hspace{0.2cm}t>0 \hspace{0.2cm} \mbox{and}\hspace{0.2cm} T_{\alpha}f(0)=\displaystyle{\lim_{t\longrightarrow0^{+}}T_{\alpha}f(t)},
\end{equation}
provided that the limits exist. In addition, if the conformable fractional derivative of $f$ of order $\alpha$ exists, then we simply say $f$ is $\alpha$-differentiable.
\end{definition}
\noindent This paper is organized to provide a comprehensive understanding of a new conformable fractional predictor-corrector algorithm and its numerical applications. In Section 2, we present the preliminaries necessary for understanding the subsequent sections, including fundamental concepts and definitions in conformable fractional calculus. Section 3 offers a detailed description of the Adams-Bashforth and Adams-Moulton methods, laying the groundwork for their extension to the fractional domain. In Section 4, we introduce the conformable fractional predictor-corrector algorithm, explaining its formulation and theoretical basis. Finally, Section 5 showcases the practical utility of the proposed algorithm through various numerical applications, demonstrating its effectiveness and accuracy in solving fractional differential equations.

\section{Preliminaries}
Recalling briefly some preliminary facts on the conformable fractional calculus.
\begin{theorem}\cite{khalil2014new} If a function $f:[0, \infty) \longrightarrow \mathbb{R}$ is $\alpha$-differentiable at $t_0>0, \alpha \in(0,1]$, then $f$ is continuous at $t_0$.
\end{theorem}
\noindent Let $\alpha\in (0,1]$, the conformable fractional derivative of certain functions \cite{khalil2014new}
\begin{multicols}{2}
\begin{itemize}
    \item $T_\alpha\left(t^p\right)=p t^{p-\alpha}$ for all $p \in \mathbb{R}$.
    \item $T_\alpha(\lambda)=0$, for all $\lambda \in \mathbb{R}$.
    \item $T_\alpha\left(e^{c t}\right)=c t^{1-\alpha} e^{c t}, c \in \mathbb{R}$.
    \item $T_\alpha(\sin b t)=b t^{1-\alpha} \cos b t, b \in \mathbb{R}$.
    \item $T_\alpha(\cos b t)=-b t^{1-\alpha} \sin b t, b \in \mathbb{R}$.
    \item $T_\alpha\left(\frac{1}{\alpha} t^\alpha\right)=1$.
    \item $T_\alpha\left(\sin \frac{1}{\alpha} t^\alpha\right)=\cos \frac{1}{\alpha} t^\alpha$.
    \item $T_\alpha\left(\cos \frac{1}{\alpha} t^\alpha\right)=-\sin \frac{1}{\alpha} t^\alpha$.
    \item $T_\alpha\left(e^{\frac{1}{\alpha} t^\alpha}\right)=e^{\frac{1}{\alpha} t^\alpha}$.
\end{itemize}
\end{multicols}
\begin{theorem}\cite{khalil2014new}
    Let $\alpha \in(0,1]$ and $f$, $g$ be $\alpha$-differentiable at a point $t>0$. Then
    \begin{enumerate}
        \item $T_\alpha(a f+b g)=a T_\alpha(f)+b T_\alpha(g)$, for all $a, b \in \mathbb{R}$.
        \item $T_\alpha(f g)=f T_\alpha(g)+g T_\alpha(f)$.
        \item $T_\alpha\left(\frac{f}{g}\right)=\frac{g T_\alpha(f)-f T_\alpha(g)}{g^2}$.
        \item If, in addition, $f$ is differentiable, then $T_\alpha(f)(t)=t^{1-\alpha} \frac{d f}{d t}(t)$.
    \end{enumerate}
\end{theorem}
\begin{definition}\cite{khalil2014new}
The conformable fractional integral $I^{\alpha}$ of a function $f$ is defined by
\begin{equation}
    I^{\alpha}f(t)=\int_{0}^{t}x^{\alpha-1}f(x)dx,\hspace{0.2cm} \mbox{for} \hspace{0.2cm}t\geq0.\label{conformable_integral}
\end{equation}
\end{definition}
\begin{theorem}\cite{khalil2014new} If $f$ is a continuous function in the domain of $I^{\alpha}$, then we have
\begin{equation}
    T_{\alpha}I^{\alpha}f(t)=f(t).
\end{equation}
\end{theorem}
\begin{theorem}\cite{abdeljawad2015conformable} If $f$ is a differentiable function, then we have
\begin{equation}
    I^{\alpha}T_{\alpha}f(t)=f(t)-f(0).
\end{equation}\label{fdiff}
\end{theorem}
\section{Description of the Adams-Bashforth and Adams-Moulton method}
\subsection{The calssical Adams-Bashforth and Adams-Moulton method}

The numerical solution of ordinary differential equations (ODEs) is a fundamental aspect of computational mathematics. Adams-Bashforth and Adams-Moulton methods are widely used for this purpose, providing efficient and accurate solutions.
\\\\
In this section, we motivate our work by recalling the classical one-step Adams-Bashforth and Adams-Moulton methods for first-order equations, also known as the predictor-corrector algorithm.
We are interested in the numerical solution of the initial value problem:

\begin{equation}
    \left \{
\begin{array}{rcl}
\frac{dy}{dt}(t)&=&f\left(t, y(t)\right)\\
y(0)&=&y_0
\end{array}
\right. \:\:\: \text{with} \:\:\:t \in [0, \tau], \quad \tau\in \mathbb{R}^{*}_{+} \label{IVP}
\end{equation}
Let $(t_j)_{j\in \{0, 1, 2, ..., n+1\}}$ be a subdivision of the interval $[0, \tau]$ with $t_0 = 0$ and $t_{n+1}= \tau$, by a simple integration over the interval $\left[t_j, t_{j+1}\right]$, we get

\begin{equation}
    \int_{t_j}^{t_{j+1}}dy(t) = \int_{t_j}^{t_{j+1}}f(t, y(t))dt.
\end{equation}
Thus one get

\begin{equation}
    y(t_{j+1}) - y(t_j) = \int_{t_j}^{t_{j+1}}f(t, y(t))dt.
\end{equation}
By letting $y_j$ denotes an approximation of $y(t_j)$ we find the following

\begin{equation}
    y_{j+1} = y_j + \int_{t_j}^{t_{j+1}}f(t, y(t))dt.\label{integralEq}
\end{equation}
The main idea behind Adams methods is approximating the integral $\int_{t_j}^{t_{j+1}}f(t, y(t))dt$ using numerical methods.\\\\
By using the two-point trapezoidal quadrature formula \cite{davis2007methods} and considering $h = (t_{j+1} - t_j)$, we obtain the following approximation

\begin{equation}
    \int_{t_j}^{t_{j+1}}f(t, y(t))dt = \frac{h}{2}\left[f(t_j, y_j) + f(t_{j+1}, y_{j+1})\right].
\end{equation}
Thus one get the following scheme
\begin{equation}
    \left \{
\begin{array}{rcl}
y_0&=&y(t_0), \\
y_{j+1}&=&y_j + \frac{h}{2}\left[f(t_j, y_j) + f(t_{j+1}, y_{j+1})\right].
\end{array}
\right. \label{AM}
\end{equation}
Formula (\ref{AM}) refers to the Adams–Moulton \cite{allahviranloo2007numerical, klopfenstein1968numerical} one-step method of order 2, note that it is an implicit scheme.\\\\
Now if we replace the trapezoidal quadrature formula by the rectangle rule one get
\begin{equation}
    \int_{t_j}^{t_{j+1}}f(t, y(t))dt = hf(t_j, y_j)
\end{equation}
Therefore the following scheme is obtained
\begin{equation}
    \left \{
\begin{array}{rcl}
y_0&=&y(t_0), \\
y_{j+1}&=&y_j + hf(t_j, y_j).
\end{array}
\right. \label{ABExplicite}
\end{equation}
Formula (\ref{ABExplicite}) refers to the Adams–Bashforth one-step method of order 2 (also known as the forward Euler method), note that it is an explicit  scheme.\\\\
Consequently if we note $y^p_{j+1}$ the predicted value of $y_{j+1}$, one can get the one-step predictor–corrector algorithm, which is a hybrid method that uses the explicit advantage of Adams–Bashforth method (\ref{ABExplicite}) and the stability \cite{hairer1991ii} of Adams–Moulton method (\ref{AM}), it is given as follows
\begin{equation}
        \left \{
\begin{array}{rcl}
y_0&=&y(t_0), \\
y^p_{j+1}&=&y_j + hf(t_j, y_j), \\
y_{j+1}&=&y_j + \frac{h}{2}\left[f(t_j, y_j) + f(t_{j+1}, y^p_{j+1})\right].
\end{array}
\right.
\end{equation}
\subsection{A predictor-corrector algorithm in the framework of Caputo type fractional derivative of order $\alpha$}
In 2002, K. Diethelm et al. in \cite{diethelm2002predictor} construct an Adams-Bashforth and Adams-Moulton method to solve the following initial value problem with Caputo type fractional derivative operator of order $\alpha$
\begin{equation}
    \left \{
\begin{array}{rcl}
^{C}D_0^\alpha y(t)&=&f\left(t, y(t)\right)\\
y^{(k)}(0)&=&y^{(k)}_0,\quad\quad k = 0, 1, ..., m-1
\end{array}
\right. \label{CFIVP}
\end{equation}
Where without loss of generality $t \in [0, \tau], \:\tau\in \mathbb{R}^{*}_{+}$, $^{C}D_0^\alpha$ is the Caputo type fractional derivative of order $\alpha>0$ defined in formula (\ref{caputoOperator}), $m = \lceil\alpha\rceil$ is just the value $\alpha$ rounded up to the nearest integer and $y^{(k)}$ is the ordinary $k$th
derivative of y. One can readily observe that the requisite number of initial conditions to determine a unique solution is denoted by $m$. A complete discussion elucidating why the authors opt to use Caputo type fractional derivative of order $\alpha$ and to use the classical ordinary derivative of $y$ as initial conditions in (\ref{CFIVP}) can be found in \cite{diethelm2002analysis, diethelm2002predictor,diethelm2004detailed}.\\\\
Note that if, (a) $f$ is continues with respect to both its arguments and (b) $f$ has a Lipschitz condition with respect to the second argument, we can indeed say that a solution exists and that
this solution is uniquely determined \cite{diethelm2002analysis}, for the initial value problem (\ref{CFIVP}).\\\\
Let $(t_j)_{j\in \{0, 1, 2, ..., n+1\}}$ be a subdivision of the interval $[0, \tau]$ with $t_0 = 0$ and $t_{n+1}= \tau$. According to \cite{diethelm2002predictor} a predictor-corrector algorithm in the framework of Caputo type fractional derivative of order $\alpha$ can be suggested as follows
\begin{equation}
        \left \{
\begin{array}{rcl}
y^{(k)}(0)&=&y^{(k)}_0\in\mathbb{R},\quad\quad k = 0, 1, ..., m-1, \\
y^p_{n+1}&=&\sum_{k=0}^{\lceil\alpha\rceil-1} \frac{t_{n+1}^k}{k !} y_0^{(k)} + \frac{h^\alpha}{\Gamma(\alpha+1)} \sum_{j=0}^{n}\left((n+1-j)^\alpha-(n-j)^\alpha\right)f(t_{j}, y_{j}), \\
y_{n+1} &=& \sum_{k=0}^{\lceil\alpha\rceil-1} \frac{t_{n+1}^k}{k !} y_0^{(k)} + \frac{h^{\alpha}}{\Gamma(\alpha+2)}\sum_{j=0}^{n} a_{j} f\left(t_j, y_{j}\right)+ \frac{h^{\alpha}}{\Gamma(\alpha+2)}f\left(t_{n+1}, y^p_{n+1}\right).
\end{array}
\right.\label{caputonalgo}
\end{equation}
Where the coefficients $a_j, \: j = 0, 1, ..., n$ are given as
\begin{equation}
    a_{j}= \begin{cases}n^{\alpha+1}-(n-\alpha)(n+1)^\alpha, & \text { if } j=0 \\ (n-j)^{\alpha+1}-2(n-j+1)^{\alpha+1}+(n-j+2)^{\alpha+1}, & \text { if } 1 \leq j \leq n \end{cases}
\end{equation}

\section{A conformable fractional Predictor-Corrector Algorithm}
Let us consider the following Initial value problem
\begin{equation}
    \left \{
\begin{array}{rcl}
T_{\alpha}y(t)&=&f\left(t, y(t)\right)\\
y(0)&=&y_0
\end{array}
\right. \:\:\: \text{with} \:\:\:t \in [0, \tau], \quad \tau\in \mathbb{R}^{*}_{+}\quad \text{and}\:\:\: \alpha \in (0, 1] \label{FIVP}
\end{equation}
Where $T_{\alpha}$  denotes the conformable fractional derivative of order $\alpha$ and $f: [0, \tau]\times\mathbb{R}\mapsto\mathbb{R}$ is a continuous function\\\\
Existence of at least one solution \cite{zhong2018basic} to the aforementioned problem was established through the application of Schauder’s fixed point theorem. Now if we suppose that $y$ is differentiable, then from theorem (\ref{fdiff}) we get
\begin{equation}
    I_{\alpha}T_{\alpha}(y)(t) = y(t) - y(0)
\end{equation}
To give an equivalent predictor-corrector algorithm in the framework of conformable fractional derivative we need to introduce a similar formula as in (\ref{integralEq}) with some unavoidable modification, in fact such formula can be obtained simply by applying the operator $I_{\alpha}$ defined in (\ref{conformable_integral}) to both sides of the initial value problem (\ref{FIVP}). Thus one obtain
\begin{equation}
    y(t) = y(0) + \int_{0}^{t}x^{\alpha-1}f(x, y(x))dx
\end{equation}
This equation exhibits a slight variation from Equation (\ref{integralEq}), as the integration range now initiates from 0 rather than $t_j$, 
however, this does not pose significant challenges in our efforts to generalize the Adams method to the conformable fractional differential equations framework. The idea is approximating the integral $\int_{0}^{t}x^{\alpha-1}g(x)dx$ using numerical methods.\\
Let $(t_j)_{j\in \{0, 1, 2, ..., n+1\}}$ be a subdivision of the interval $[0, \tau]$ with $t_0 = 0$ and $t_{n+1}= \tau$.\\\\
One can  simply use the product trapezoidal quadrature formula to approximate the integral
\begin{equation}
    \int_0^{t_{n+1}}x^{\alpha-1} g(x) \mathrm{d} x \approx \int_0^{t_{n+1}}x^{\alpha-1} \Tilde{g}(x) \mathrm{d} z,\label{integralForApproximate}
\end{equation}
Here, $\tilde{g}$ represents the piecewise linear interpolant for $g$, with nodes and knots selected at $t_j$, $j=0,1,2, \ldots, n+1$. Consequently
\begin{equation}
    \int_0^{t_{n+1}}x^{\alpha-1} g(x) \mathrm{d} x \approx \frac{h^{\alpha}}{\alpha(\alpha+1)}\sum_{j=0}^{n+1} a_{j} g\left(t_j\right),\label{FAM1}
\end{equation}
where
\begin{equation}
    a_{j}= \begin{cases}1, & \text { if } j=0 \\ (j-1)^{\alpha+1}-2j^{\alpha+1}+(j+1)^{\alpha+1}, & \text { if } 1 \leq j \leq n \\ (\alpha+1)(n+1)^{\alpha}+n^{\alpha+1}-(n+1)^{\alpha+1}, & \text { if } j=n+1\end{cases}\label{mainaj}
\end{equation}
Hence if we note $y_j$ an approximation of $y(t_j)$ for $t_j$, $j=0,1,2, \ldots, n+1$ one get the following scheme
\begin{equation}
\begin{aligned}
    y_{n+1} &= y_0 + \frac{h^{\alpha}}{\alpha(\alpha+1)}\sum_{j=0}^{n} a_{j} f\left(t_j, y_{j}\right) \\
    &\quad + \frac{h^{\alpha}}{\alpha(\alpha+1)}\left((\alpha+1)(n+1)^{\alpha}+n^{\alpha+1}-(n+1)^{\alpha+1}\right)f\left(t_{n+1}, y_{n+1}\right),
\end{aligned} \label{FAM2}
\end{equation}
As a result formula (\ref{FAM2}) refers to the conformable fractional Adams–Moulton method, note that it is an implicit scheme, which means we need to predict the value of $y_{n+1}$ an approximation of $y(t_{n+1})$, to do that it is a must to introduce a conformable fractional Adams–Bashforth method.\\\\
By using the product rectangle rule we can approximate the integral on the left-hand side of Equation (\ref{integralForApproximate}) as
\begin{equation}
    \int_0^{t_{n+1}}x^{\alpha-1} g(x) \mathrm{d} x \approx \sum_{j=0}^{n} \int_{t_{j}}^{t_{j+1}}x^{\alpha-1} g(t_{j}) \mathrm{d} x = \frac{h^\alpha}{\alpha} \sum_{j=0}^{n}\left((j+1)^\alpha - j^\alpha\right)g(t_{j}).
\end{equation}
Consequently if we note $y^p_{n+1}$ the predicted value of $y_{n+1}$, one get the conformable fractional Adams–Bashforth method
\begin{equation}
    y^p_{n+1} = y_0 + \frac{h^\alpha}{\alpha} \sum_{j=0}^{n}\left((j+1)^\alpha - j^\alpha\right)f(t_{j}, y_{j}).
\end{equation}
As a result one get a conformable fractional predictor-corrector algorithm
\begin{equation}
        \left \{
\begin{array}{rcl}
y_0&=&y(0), \\
y^p_{n+1}&=&y_0 + \frac{h^\alpha}{\alpha} \sum_{j=0}^{n}\left((j+1)^\alpha - j^\alpha\right)f(t_{j}, y_{j}), \\
y_{n+1} &=& y_0 + \frac{h^{\alpha}}{\alpha(\alpha+1)}\sum_{j=0}^{n} a_{j} f\left(t_j, y_{j}\right)+ \frac{h^{\alpha}}{\alpha(\alpha+1)}\left((\alpha+1)(n+1)^{\alpha}+n^{\alpha+1}-(n+1)^{\alpha+1}\right)f\left(t_{n+1}, y^p_{n+1}\right).
\end{array}
\right.\label{themainalgo}
\end{equation}
Where the coefficients $a_j, \: j = 0, 1, ..., n$ are given as
\begin{equation}
    a_{j}= \begin{cases}1, & \text { if } j=0 \\ (j-1)^{\alpha+1}-2j^{\alpha+1}+(j+1)^{\alpha+1}, & \text { if } 1 \leq j \leq n \end{cases}
\end{equation}
\section{Numerical application}
\begin{example}
    To give a concrete example of the application of our algorithm, one can consider the following conformable fractional initial value problem 
\begin{equation}
    \left \{
\begin{array}{rcl}
T_{\alpha}y(t)&=&ty(t)\\
y(0)&=&1
\end{array}
\right. \:\:\: \text{with} \:\:\:t \in [0, \tau], \quad \tau\in \mathbb{R}^{*}_{+}\quad \text{and}\:\:\: \alpha \in (0, 1] \label{FIVPexample1}
\end{equation}
\begin{lemma}
    Let $\alpha \in (0, 1]$ and $y$ be $\alpha$-differentiable  at a point $t>0$, then if, in addition, $y$ is differentiable \cite{khalil2014new}, one get 
    \begin{equation}
        T_{\alpha}y(t) = t^{1-\alpha}\frac{dy}{dt}(t).
    \end{equation}
    \label{lemma1}
\end{lemma}
\noindent One can easily use lemma (\ref{lemma1}), to verify that the exact solution of this initial value problem, has the form
\begin{equation}
    y(t) = \exp(\frac{t^{\alpha+1}}{\alpha+1}).
\end{equation}
By using the algorithm (\ref{themainalgo}), one can have the plot in (Figure \ref{fig:example1}) (obtained using the code in Appendix A), of the numerical and exact solution, of the conformable fractional initial value problem (\ref{FIVPexample1}).
\end{example}

\begin{figure}[H]
    \centering
    \includegraphics[scale=.8]{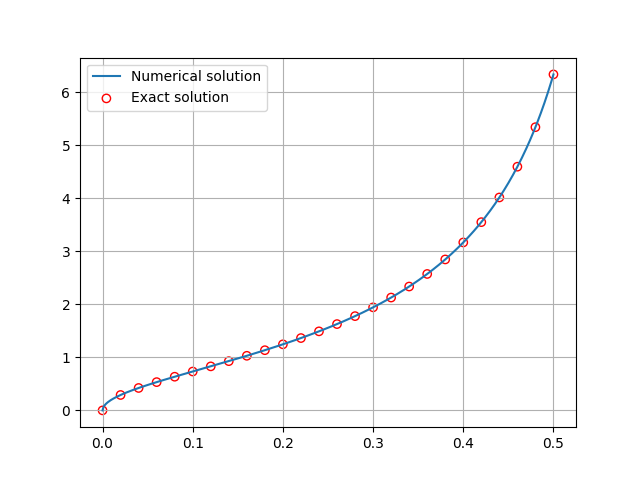} %
    \caption{The plot of $y(t)$, the solution of the conformable fractional initial value problem (\ref{FIVPexample2}), for $\alpha=0.5$.}
    \label{fig:example2}
\end{figure}

\begin{example}
    Consider the following non-linear conformable fractional differential equation
    \begin{equation}
    \left \{
\begin{array}{rcl}
T_{\alpha}y(t)&=&1 + y^2(t)\\
y(0)&=&0
\end{array}
\right. \:\:\: \text{with} \:\:\:t \in [0, \tau], \quad \tau\in \mathbb{R}^{*}_{+}\quad \text{and}\:\:\: \alpha \in (0, 1] \label{FIVPexample2}
\end{equation}
Recalling lemma (\ref{lemma1}), separating variables and integrating, one get $\arctan(y(t)) = \frac{t^\alpha}{\alpha}$. Using the fact that the range of the inverse tangent function is $(-\frac{\pi}{2}, \frac{\pi}{2})$, we get $0 \leq \frac{t^\alpha}{\alpha} < \frac{\pi}{2}$. Therefore the domain of the solution to the initial value problem (\ref{FIVPexample2}) is the interval $0 \leq t < \left(\alpha\frac{\pi}{2}\right)^{\frac{1}{\alpha}}$
\\\\
For this reason it follows that the exact solution of (\ref{FIVPexample2}) is 
\begin{equation}
    y(t) = \tan(\frac{t^\alpha}{\alpha})\:,\quad\:\: \text{with}\quad t \in \left[0,\: \left(\alpha\frac{\pi}{2}\right)^{\frac{1}{\alpha}}\right)
\end{equation}
Note that the solution cannot
jump over the vertical asymptote. To avoid the singularity point of the solution $y$, one can plot the numerical and exact solution, of the conformable fractional initial value problem (\ref{FIVPexample2}) over the interval $[0, \frac{1}{2}]$ for $\alpha=0.5$, note that here $\tau = \frac{1}{2}$, See (Figure \ref{fig:example2}).
\end{example}

\begin{figure}[H]
    \centering
    \includegraphics[scale=.8]{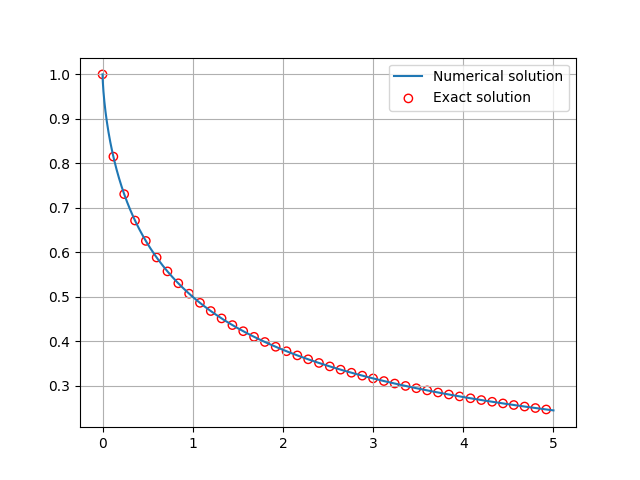} %
    \caption{The plot of $y(t)$, the solution of the conformable fractional initial value problem (\ref{FIVPexample3}), for $\alpha=0.7$.}
    \label{fig:example3}
\end{figure}

\begin{example}
    It is worth providing another example to demonstrate the proficiency of our algorithm in generating numerical solutions that closely approximate the exact solution, with negligible error that is imperceptible.\\\\
    One can consider the following non-linear conformable fractional differential equation
    \begin{equation}
    \left \{
\begin{array}{rcl}
T_{\alpha}y(t)&=&-\alpha y^2(t)\\
y(0)&=&1
\end{array}
\right. \:\:\: \text{with} \:\:\:t \in [0, \tau], \quad \tau\in \mathbb{R}^{*}_{+}\quad \text{and}\:\:\: \alpha \in (0, 1] \label{FIVPexample3}
\end{equation}
One can easily verify that the exact solution of the conformable fractional differential equation (\ref{FIVPexample3}) can be written as
\begin{equation}
    y(t) = \frac{1}{1 + t^\alpha}
\end{equation}
The plot of the numerical solution and the exact solution of the conformable fractional differential equation (\ref{FIVPexample3}) is given in (Figure \ref{fig:example3}), for as an example $\alpha = 0.7$.
\end{example}
\begin{figure}[H]
    \centering
    \includegraphics[scale=.8]{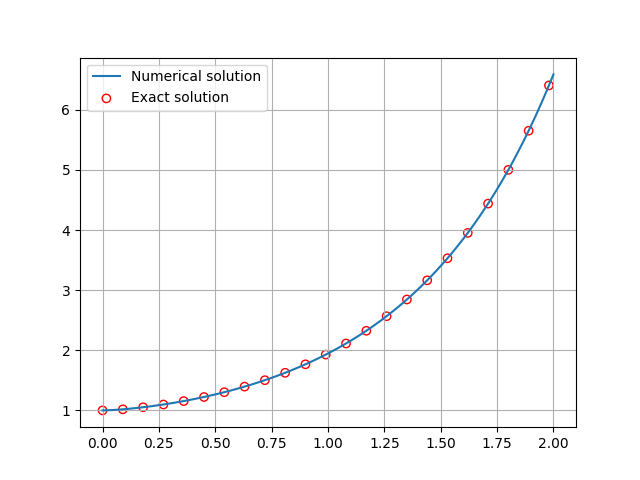} %
    \caption{The plot of $y(t)$, the solution of the conformable fractional initial value problem (\ref{FIVPexample1}), for $\alpha=0.5$.}
    \label{fig:example1}
\end{figure}

\newpage
{\Large\textbf{Appendix A: A Python Code of the Algorithm (For Example (\ref{FIVPexample1}))}}
\\\\
\begin{lstlisting}
from numpy import *
from matplotlib.pyplot import *

def ConformableAm(f, a, b, y0, h, alpha):
    n = int((b - a) / h)
    T = [a + i * h for i in range(n + 1)]
    Y = [0] * (n + 1)
    cte1 = pow(h, alpha) / alpha
    cte2 = cte1 / (alpha + 1)
    Y[0] = y0
    k = Y[0] + cte1 * f(T[0], Y[0])
    Y[1] = Y[0] + cte2 * (f(T[0], Y[0]) + alpha * f(T[1], k))

    s = Y[0] + cte2 * f(T[0], Y[0])
    for i in range(1, n):
        k += cte1 * ((i+1)**alpha - i**alpha) * f(T[i], Y[i])
        s += ((i-1)**(alpha+1) - 2 * (i**(alpha+1)) 
                + (i+1)**(alpha+1)) * f(T[i], Y[i]) * cte2
        Y[i+1] = s + cte2 * ((alpha+1)*((i+1)**alpha) + i**(alpha+1) 
                                - (i+1)**(alpha+1)) * f(T[i+1], k)

    return T, Y

def f(t, y):
    return t * y

alpha = 0.5

T, Y = ConformableAm(f, 0, 2, 1, 0.001, alpha)

Z = [0] * len(T)

for i in range(len(Z)):
    Z[i] = exp((T[i] ** (alpha + 1)) / (alpha + 1))

plot(T, Y, label="Numerical solution")
# Plot Exact Solution with Markers
tt = [T[i] for i in range(0, len(T), 90)]
yy = [Z[i] for i in range(0, len(T), 90)]
scatter(tt, yy, marker='o', facecolor='none', edgecolor='red',
        label="Exact solution")

legend()
grid()
show()
\end{lstlisting}

\newpage
\bibliographystyle{unsrt} 
\bibliography{Mybib.bib}

\begin{thebibliography}{10}

\bibitem{butzer2000introduction}
Paul~L Butzer and Ursula Westphal.
\newblock An introduction to fractional calculus.
\newblock In {\em Applications of fractional calculus in physics}, pages 1--85.
  World Scientific, 2000.

\bibitem{machado2011recent}
J~Tenreiro Machado, Virginia Kiryakova, and Francesco Mainardi.
\newblock Recent history of fractional calculus.
\newblock {\em Communications in nonlinear science and numerical simulation},
  16(3):1140--1153, 2011.

\bibitem{shaikh2020mathematical}
Amjad~Salim Shaikh, Iqbal~Najiroddin Shaikh, and Kottakkaran~Sooppy Nisar.
\newblock A mathematical model of covid-19 using fractional derivative:
  outbreak in india with dynamics of transmission and control.
\newblock {\em Advances in Difference Equations}, 2020(1):373, 2020.

\bibitem{magin2010fractional}
Richard~L Magin.
\newblock Fractional calculus models of complex dynamics in biological tissues.
\newblock {\em Computers \& Mathematics with Applications}, 59(5):1586--1593,
  2010.

\bibitem{dalir2010applications}
Mehdi Dalir and Majid Bashour.
\newblock Applications of fractional calculus.
\newblock {\em Applied Mathematical Sciences}, 4(21):1021--1032, 2010.

\bibitem{koh1990application}
Chan~Ghee Koh and James~M Kelly.
\newblock Application of fractional derivatives to seismic analysis of
  base-isolated models.
\newblock {\em Earthquake engineering \& structural dynamics}, 19(2):229--241,
  1990.

\bibitem{assaleh2007modeling}
Khaled Assaleh and Wajdi~M Ahmad.
\newblock Modeling of speech signals using fractional calculus.
\newblock In {\em 2007 9th International Symposium on Signal Processing and Its
  Applications}, pages 1--4. IEEE, 2007.

\bibitem{mainardi2022fractional}
Francesco Mainardi.
\newblock {\em Fractional calculus and waves in linear viscoelasticity: an
  introduction to mathematical models}.
\newblock World Scientific, 2022.

\bibitem{bagley1983theoretical}
Ronald~L Bagley and Peter~J Torvik.
\newblock A theoretical basis for the application of fractional calculus to
  viscoelasticity.
\newblock {\em Journal of rheology}, 27(3):201--210, 1983.

\bibitem{meral2010fractional}
FC~Meral, TJ~Royston, and R~Magin.
\newblock Fractional calculus in viscoelasticity: an experimental study.
\newblock {\em Communications in nonlinear science and numerical simulation},
  15(4):939--945, 2010.

\bibitem{koeller1984applications}
RC747787 Koeller.
\newblock Applications of fractional calculus to the theory of viscoelasticity.
\newblock 1984.

\bibitem{fellah2002application}
Zine El~Abiddine Fellah, C~Depollier, and Mohamed Fellah.
\newblock Application of fractional calculus to the sound waves propagation in
  rigid porous materials: validation via ultrasonic measurements.
\newblock {\em Acta Acustica united with Acustica}, 88(1):34--39, 2002.

\bibitem{mathieu2003fractional}
Beno{\^\i}t Mathieu, Pierre Melchior, Alain Oustaloup, and Ch~Ceyral.
\newblock Fractional differentiation for edge detection.
\newblock {\em Signal Processing}, 83(11):2421--2432, 2003.

\bibitem{kulish2002application}
Vladimir~V Kulish and Jos{\'e}~L Lage.
\newblock Application of fractional calculus to fluid mechanics.
\newblock {\em J. Fluids Eng.}, 124(3):803--806, 2002.

\bibitem{podlubny1998fractional}
Igor Podlubny.
\newblock {\em Fractional differential equations: an introduction to fractional
  derivatives, fractional differential equations, to methods of their solution
  and some of their applications}.
\newblock elsevier, 1998.

\bibitem{khalil2014new}
Roshdi Khalil, Mohammed Al~Horani, Abdelrahman Yousef, and Mohammad Sababheh.
\newblock A new definition of fractional derivative.
\newblock {\em Journal of computational and applied mathematics}, 264:65--70,
  2014.

\bibitem{abdeljawad2015conformable}
Thabet Abdeljawad.
\newblock On conformable fractional calculus.
\newblock {\em Journal of computational and Applied Mathematics}, 279:57--66,
  2015.

\bibitem{davis2007methods}
Philip~J Davis and Philip Rabinowitz.
\newblock {\em Methods of numerical integration}.
\newblock Courier Corporation, 2007.

\bibitem{allahviranloo2007numerical}
Tofigh Allahviranloo, Nazanin Ahmady, and E~Ahmady.
\newblock Numerical solution of fuzzy differential equations by
  predictor--corrector method.
\newblock {\em Information sciences}, 177(7):1633--1647, 2007.

\bibitem{klopfenstein1968numerical}
RW~Klopfenstein and RS~Millman.
\newblock Numerical stability of a one-evaluation predictor-corrector algorithm
  for numerical solution of ordinary differential equations.
\newblock {\em Mathematics of Computation}, 22(103):557--564, 1968.

\bibitem{hairer1991ii}
E~Hairer, G~Wanner, and O~Solving.
\newblock Ii: Stiff and differential-algebraic problems.
\newblock {\em Berlin [etc.]: Springer}, 1991.

\bibitem{diethelm2002predictor}
Kai Diethelm, Neville~J Ford, and Alan~D Freed.
\newblock A predictor-corrector approach for the numerical solution of
  fractional differential equations.
\newblock {\em Nonlinear Dynamics}, 29:3--22, 2002.

\bibitem{diethelm2002analysis}
Kai Diethelm and Neville~J Ford.
\newblock Analysis of fractional differential equations.
\newblock {\em Journal of Mathematical Analysis and Applications},
  265(2):229--248, 2002.

\bibitem{diethelm2004detailed}
Kai Diethelm, Neville~J Ford, and Alan~D Freed.
\newblock Detailed error analysis for a fractional adams method.
\newblock {\em Numerical algorithms}, 36:31--52, 2004.

\bibitem{zhong2018basic}
Wenyong Zhong and Lanfang Wang.
\newblock Basic theory of initial value problems of conformable fractional
  differential equations.
\newblock {\em Advances in Difference Equations}, 2018:1--14, 2018.

\end{thebibliography}

\end{document}